\documentclass[reqno]{amsart}
\usepackage{amsmath}
\allowdisplaybreaks[3]
\numberwithin{equation}{section}
\numberwithin{equation}{section}

\def\proof{\indent{\em Proof.\quad}}
\def\endproof{\hfill\hbox{$\sqcup$}\llap{\hbox{$\sqcap$}}\medskip}
\newtheorem{thm}{{\indent\bf Theorem}}[section]
\newtheorem{prop}{{\indent\bf Proposition}}[section]
\newtheorem{dfn}{{\indent\bf Definition}}[section]
\newtheorem{rmk}{{\indent\bf Remark}}[section]
\newtheorem{lem}{{\indent\bf Lemma}}[section]
\newtheorem{cor}{{\indent\bf Corollary}}[section]
\newtheorem{expl}{{\indent\bf Example}}[section]

\newcommand{\mb}{\mbox}
\newcommand{\hs}{\hspace}

\newcommand{\ol}{\overline}

\newcommand{\strl}{\stackrel}

\newcommand{\td}{\tilde}

\newcommand{\fr}{\frac}
\newcommand{\ed}{{\rm End}}
\newcommand{\edd}{\end{document}}
\newcommand{\be}{\begin{equation}}
\newcommand{\ee}{\end{equation}}

\newcommand{\lagl}{\langle}
\newcommand{\ragl}{\rangle}
\newcommand{\lmx}{\left(\begin{matrix}}
\newcommand{\rmx}{\end{matrix}\right)}
\newcommand{\ldt}{\left|\begin{matrix}}
\newcommand{\rdt}{\end{matrix}\right|}

\newcommand{\rank}{{\rm rank\,}}

\newcommand{\tr}{{\rm tr\,}}
\newcommand{\vfi}{\varphi}

\newcommand{\bbr}{{\mathbb R}}

\newcommand{\bbc}{{\mathbb C}}

\newcommand{\ba}{\begin{array}}
\newcommand{\ea}{\end{array}}
\newcommand{\nnm}{\nonumber}
\newcommand{\beal}{\begin{align}}
\newcommand{\eal}{\end{align}}
\newcommand{\bea}{\begin{eqnarray}}
\newcommand{\eea}{\end{eqnarray}}

\newcommand{\pp}[2]{\fr{\partial #1}{\partial #2}}

\title[A new characterization of Calabi composition]
{A new characterization of Calabi composition\\of hyperbolic affine hyperspheres}%
\author{Xingxiao Li}%

\begin{document}

\begin{abstract}
In this paper, we mainly prove a theorem with a corollary establishing two characterizations of the Calabi composition of hyperbolic hyperspheres, where the second characterization (i.e., the corollary) has been given via a dual correspondence theorem earlier but now we would like to use a very direct method. Note that Z.J. Hu, H.Z. Li and L. Vrancken also gave a characterization of the $2$-factor Calabi composition in a different manner.

\noindent
keywords and expressions: hyperbolic affine spheres;Calabi composition;characterization;reducibility of Riemannian manifolds

\noindent
AMS Classification: 53C
\end{abstract}
\maketitle
% ----------------------------------------------------------------

\tableofcontents

\section{Introduction}

As we know, affine hyperspheres are the most important objects studied in affine differential geometry of hypersursurfaces, drawing great attention of many geometers.
In fact, affine hyperspheres seems simple in definition but they do form a very large class of hypersurfaces, the study of which is fruitful in recent twenty years. For example, the proof of the Calabi's conjecture (\cite{amli90}, \cite{amli92}), the classification of hyperspheres of constant affine curvatures (\cite{vra-li-sim91},  \cite{wang93}, \cite{kri-vra99}), and in \cite{hu-li-vra11} the complete classification of locally strongly convex hypersurfaces with parallel Fubini-Pick forms as a special class of hyperbolic affine hyperspheres (for some earlier partial results, see \cite{bok-nom-sim90}, \cite{dil-vra-yap94}, \cite{hu-li-sim-vra09}). As for the general nondegenerate case, there also have been some interesting partial classification results, see for example the series of published papers by Z.J. Hu etc: \cite{hu-li11}, \cite{hu-li-li-vra11a} and \cite{hu-li-li-vra11b}.

In 1972, E. Calabi (\cite{cal72}) found a composition formula by which one can construct new hyperbolic affine hyperspheres from any two given ones. The present author has generalized Calabi construction to the case of multiple factors (See \cite{lix93}, published in Chinese). Later in 1994 F. Dillen and L. Vrancken \cite{dil-vra94} generalized Calabi original composition to any two proper affine hyperspheres and gave a detailed study of these composed affine hyperspheres. They also mentioned that their construction applies to the case of multiple factors but with no detail of this. In 2008, in order to establish their later classification in \cite{hu-li-vra11} mentioned above, Z.J. Hu, H.Z. Li and L. Vrancken proved a characterization of the Calabi composition of hyperbolic hyperspheres (\cite{hu-li-vra08}) by special decompositions of the tangent bundle. We would like to remark that, by using the similar idea of \cite{hu-li-vra08}, H.Z. Li and X.F. Wang has in a way characterized the so called Calabi product of parallel Lagrangian submanifolds in the complex projective space $\bbc P^n$ (\cite{li-wang11}).

In a previous paper, we explicitly defined the Calabi composition of multiple factors of hyperbolic hyperspheres, possibly including some point factors viewing as ``$0$-dimensional hyperbolic hyperspheres'', and made it in detail for the computation of the basic affine invariants of this composition. In this article, by using those basic affine invariants,  we prove a theorem (see Theorem \ref{the main}) which provides a new and more natural characterization of the Calabi composition of multiple hyperbolic hyperspheres. In the case of affine symmetric factors this characterization turns out to be much more simple (see Corollary \ref{cor}).

{\sc Acknowledgement} The author is grateful to Professor A-M Li for his encouragement and important suggestions during the preparation of this article. He also thanks Professor Z.J. Hu for providing him valuable related references some of which are listed in the end of this paper.

\section{Preliminaries}

\subsection{The equiaffine geometry of hypersurfaces}

In this subsection, we brief some basic facts in the equiaffine geometry of hypersurfaces. For details the readers are referred to some text books, say, \cite{li-sim-zhao93} and \cite{nom-sas94}.

Let $x:M^n\to\bbr^{n+1}$ be a nondegenerate hypersurface. Then there are several basic equiaffine invariants of $x$ among which are: the affine metric (Berwald-Blaschke metric) $g$, the affine normal $\xi:=\fr1n\Delta_gx$, the Fubini-Pick $3$-form (the so called cubic form) $A\in\bigodot^3T^*M^n$ and the affine second fundamental $2$-form $B\in\bigodot^2T^*M^n$. By using the index lifting by the metric $g$, we can identify $A$ and $B$ with the linear maps $A:TM^n\to \ed(TM^n)$ or $A:TM^n\bigodot TM^n\to TM^n$ and $B:TM^n\to TM^n$, respectively, by
\be\label{ab}
g(A(X)Y,Z)=A(X,Y,Z) \mb{\ or\ }g(A(X,Y),Z)=A(X,Y,Z),\quad
g(B(X),Y)=B(X,Y),
\ee
for all $X,Y,Z\in TM^n$. Sometimes we call the corresponding $B\in \ed(TM^n)$ the affine shape operator of $x$. In this sense, the affine Gauss equation can be written as follows:
\be\label{gaus}
R(X,Y)Z=\fr12(g(Y,Z)B(X)+B(Y,Z)X-g(X,Z)B(Y)-B(X,Z)Y)-[A(X),A(Y)](Z),
\ee
where, for any linear transformations $T,S\in \ed(TM^n)$,
\be\label{comm}
[T,S]=T\circ S-S\circ T.
\ee
Each of the eigenvalues $B_1,\cdots,B_n$ of the linear map $B:TM^n\to TM^n$ is called the affine principal curvature of $x$. Define
\be\label{afme}
L_1:=\fr1n\tr B=\fr1n\sum_iB_i.
\ee
Then $L_1$ is referred to as the affine mean curvature of $x$. A hypersurface $x$ is called an (elliptic, parabolic, or hyperbolic) affine hypersphere, if all of its affine principal curvatures are equal to one (positive, 0, or negative) constant. In this case we have
\be\label{afsp}
B(X)=L_1X,\quad\mb{for all\ }X\in TM^n.
\ee
It follows that the affine Gauss equation \eqref{gaus} of an affine hypersphere assumes the following form:
\be\label{gaus_af sph}
R(X,Y)Z=L_1(g(Y,Z)X-g(X,Z)Y)-[A(X),A(Y)](Z),
\ee

Furthermore, all the affine lines of an elliptic affine hypersphere or a hyperbolic affine hypersphere $x:M^n\to\bbr^{n+1}$ pass through a fix point $o$ which is refer to as the affine center of $x$; Both the elliptic affine hyperspheres and the hyperbolic affine hyperspheres are called proper affine hyperspheres, while the parabolic affine hyperspheres are called improper affine hyperspheres.

For each vector field $\eta$ transversal to the tangent space of $x$, we have the following direct decomposition
$$
x^*T\bbr^{n+1}=x_*(TM)\oplus \bbr\cdot\eta.
$$
This decomposition and the canonical differentiation $\bar D^0$ on $\bbr^{n+1}$ define a bilinear form $h\in\bigodot^2T^*M^n$ and a connection $D^\eta$ on $TM^n$ as follows:
\be\label{dfn h}
\bar D^0_XY=x_*(D^\eta_XY)+h(X,Y)\eta,\quad\forall X,Y\in TM^n.
\ee
\eqref{dfn h} can be referred as to the affine Gauss formula of the hypersurface $x$.
In particular, in case that $\eta$ is parallel to the affine normal $\xi$, the induced connection $\nabla:=D^\eta$ is independent of the choice of $\eta$ and is referred to as the affine connection of $x$.

In what follows we make the following convention for the range of indices:
$$1\leq i,j,k,l\leq n.$$

Let $\{e_i,e_{n+1}\}$ be a local unimodular frame field along $x$ with $\eta:=e_{n+1}$ parallel to the affine normal $\xi$, and $\{\omega^i,\omega^{n+1}\}$ be its dual coframe. Then we have connection forms $\omega^A_B$, $1\leq A,B\leq n+1$, defined by
$$
d\omega^A=\omega^B\wedge\omega^A_B,\quad d\omega^A_B=\sum_{C=1}^{n+1}\omega^C_A\wedge\omega^B_C,\quad \omega^{n+1}\equiv 0.
$$
Furthermore, the above invariants can be respectively expressed locally as
\be\label{gab}
g=\sum g_{ij}\omega^i\omega^j,\quad A=\sum A_{ijk}\omega^i\omega^j\omega^k,\quad B=\sum B_{ij}\omega^i\omega^j,
\ee
subject to the following basic formulas:
\begin{align}
&\sum_{i,j} g^{ij}A_{ijk}=0\text{\ (the apolarity)},\label{basic1}\\
&A_{ijk,l}-A_{ijl,k}=\fr12(g_{ik}B_{jl}+g_{jl}B_{ik} -g_{il}B_{jk}-g_{jk}B_{il}),\label{basic3}\\
&\sum_{l}A^l_{ij,l}=\fr n2(L_1g_{ij}-B_{ij}),\label{basic3-1}
\end{align}
where $A_{ijk,l}$ are the covariant derivatives of $A_{ijk}$ with respect to the Levi-Civita connection of $g$.

Write $h=\sum h_{ij}\omega^i\omega^j$ and $H=|\det(h_{ij})|$. Then
\be\label{dfn g}
g_{ij}=H^{-\fr1{n+2}}h_{ij},\quad \xi=H^{\fr1{n+2}}e_{n+1}.
\ee
Define
\be\label{hijk0}
\sum_kh_{ijk}\omega^k=dh_{ij}+h_{ij}\omega^{n+1}_{n+1}-\sum h_{kj}\omega^k_i-\sum h_{ik}\omega^k_j.
\ee
Then the Fubini-Pick form $A$ can be determined by the following formula:
\be\label{hijktoaijk}
A_{ijk}=-\fr12H^{-\fr1{n+2}}h_{ijk}.
\ee

Define the normalized scalar curvature $\chi$ and the Pick invariant $J$ by
$$
\chi=\fr1{n(n-1)}\sum g^{il}g^{jk}R_{ijkl},\quad J=\fr1{n(n-1)}\sum A_{ijk}A_{pqr}g^{ip}g^{jq}g^{kr}.$$
Then the affine Gauss equation can be written in terms of the metric and the Fubini-Pick form as follows
\begin{align}
R_{ijkl}=&(A_{ijk,l}-A_{ijl,k})+(\chi-J)(g_{il}g_{jk}-g_{ik}g_{jl})\nnm\\ &\ +\fr2n\sum(g_{ik}A_{jlm,m}-g_{il}A_{jkm,m}) +\sum_m(A^m_{ik}A_{jlm}-A^m_{il}A_{jkm}).
\label{basic2}\end{align}

The following existence and uniqueness theorems are well known:

{\thm\label{affine existence} $($\cite{li-sim-zhao93}$)$
$($The existence$)$ Let $(M^n,g)$ be a simply connected Riemannian manifold
of dimension $n$, and $A$ be a symmetric $3$-form on $M^n$ satisfying the
affine Gauss equation \eqref{basic2} and the apolarity condition \eqref{basic1}. Then there exists a locally strongly convex immersion $x:M^n\to \bbr^{n+1}$ such that $g$ and $A$ are the affine metric and the Fubini-Pick form for $x$, respectively.}

{\thm\label{affine uniqueness} $($\cite{li-sim-zhao93}$)$ $($The uniqueness$)$ Let $x:M^n\to \bbr^{n+1}$,
$\bar x:\bar M^n\to \bbr^{n+1}$ be two locally strongly convex hypersurfaces of dimension $n$ with respectively the affine metrics $g$, $\bar g$ and the Fubini-Pick forms $A$, $\bar A$, and $\vfi:(M^n,g)\to (\bar M^n,\bar g)$ be an isometry between Riemannian manifolds. Then $\vfi^*\bar A=A$ if and only if there exists a unimodular affine transformation $\Phi:\bbr^{n+1}\to \bbr^{n+1}$ such that $\bar x\circ\vfi=\Phi\circ x$, or equivalently, $\bar x=\Phi\circ x\circ\vfi^{-1}$.}

\rmk\rm The necessity part of Theorem \ref{affine uniqueness} is proved in \cite{li-sim-zhao93}. Here we give a proof for the sufficient part as follows:

Choose an orthonormal frame field $\{e_i;\ 1\leq i\leq n\}$ on $M^n$ with its dual coframe $\{\omega^i;\ 1\leq i\leq n\}$. Let $\xi,\bar \xi$ be respectively the affine normal of $x$ and $\bar x$. Then $\{e_1,\cdots,e_n,\xi\}$ is unimodular.
Define $\bar e_i=\vfi_*(e_i)$, $\bar\omega^i=(\vfi^{-1})^*\omega^i$, $1\leq i\leq n$. Then $\{\bar\omega^i;\ 1\leq i\leq n\}$ is the dual coframe of $\{\bar e_i;\ 1\leq i\leq n\}$. Since $\vfi$ is an isometry, $\{\bar e_1,\cdots,\bar e_n,\bar \xi\}$ is also unimodular.

Under the condition that $\bar x=\Phi\circ x\circ\vfi^{-1}$, we claim that $\bar \xi=(\Phi_*(\xi))\circ {\vfi^{-1}}$. In fact
\be\label{bar ejei}
\bar e_j(\bar e_i\bar x)=\vfi_*(e_j)(\vfi_*(e_i)(\Phi\circ x\circ\vfi^{-1})) =\vfi_*(e_j)
((e_i(\Phi\circ x))\circ {\vfi^{-1}})
=(e_j(e_i(\Phi\circ x)))\circ {\vfi^{-1}}.
\ee
Denote respectively by $\nabla,\hat\nabla,\Delta$ and $\bar\nabla,\hat{\bar\nabla},\bar\Delta$ the affine connections of $x,\bar x$, the Riemannian connections and the Laplacians of $g,\bar g$. Then we find
\begin{align}
\bar \xi=&\fr1n\bar\Delta\bar x=\fr1n\left(\sum_i\left(\bar e_i(\bar e_i\bar x)-(\hat{\bar\nabla}_{\bar e_i}\bar e_i)(x)\right)\right)\nnm\\
=&\fr1n\left(\sum_i\left((e_i(e_i(\Phi\circ x)))\circ {\vfi^{-1}}-\vfi_*(\hat\nabla_{e_i}e_i)(\bar x)\right)\right)\nnm\\
=&\fr1n\left(\sum_i\left((e_i(\Phi_*(e_ix)))\circ {\vfi^{-1}} -(\hat\nabla_{e_i}e_i)(\Phi\circ x)\circ{\vfi^{-1}}\right)\right)\nnm\\
=&\fr1n\left(\sum_i\left((\Phi_*(e_i(e_ix)))\circ {\vfi^{-1}} -(\Phi_*(\hat\nabla_{e_i}e_i)(x))\circ {\vfi^{-1}}\right)\right)\nnm\\
=&\fr1n\Phi_*\left(\sum_i\left(e_i(e_ix) -(\hat\nabla_{e_i}e_i)(x)\right)\right)\circ {\vfi^{-1}} \nnm\\
=&\fr1n(\Phi_*(\Delta x))\circ {\vfi^{-1}}
=(\Phi_*(\xi))\circ {\vfi^{-1}}.\nnm
\end{align}

On the other hand, by \eqref{bar ejei} and the affine Gauss formula \eqref{dfn h} of $x$
\begin{align}
\bar e_j\bar e_i \bar x=&(e_je_i(\Phi\circ x))\circ\vfi^{-1} =(e_j\Phi_*(e_i(x))\circ\vfi^{-1} =(\Phi_*(e_je_i(x))\circ\vfi^{-1}\nnm\\
=&\left(\Phi_*\left(x_*(\nabla_{e_j}e_i)+\delta_{ij}\xi\right)\right)\circ\vfi^{-1}\nnm\\
=&\left(\Phi_*\left(x_*(\nabla_{e_j}e_i)\right)\right)\circ\vfi^{-1} +\delta_{ij}\left(\Phi_*(\xi)\right)\circ\vfi^{-1}.\nnm
\end{align}
But, by the affine Gauss formula of $\bar x$,
$$
\bar e_j\bar e_i \bar x=\bar x_*(\bar\nabla_{\bar e_j}\bar e_i)+\delta_{ij}\bar \xi =\Phi_*\left(x_*(\vfi^{-1}_*(\bar\nabla_{\bar e_j}\bar e_i))\right)\circ\vfi^{-1} +\delta_{ij}(\Phi_*(\xi))\circ\vfi^{-1}.
$$
It follows that
$$
\Phi_*\left(x_*(\nabla_{e_j}e_i)\right)=\Phi_*\left(x_*(\vfi^{-1}_*(\bar\nabla_{\bar e_j}\bar e_i))\right).
$$
Therefore
\be
\vfi^{-1}_*(\bar\nabla_{\bar e_j}\bar e_i)=\nabla_{e_j}e_i,
\text{\ or equivalently,\ }
\vfi_*(\nabla_{e_j}e_i)=\bar\nabla_{\bar e_j}\bar e_i,
\ee
from which we find that
\begin{align}
\bar A(\bar e_i,\bar e_j,\bar e_k)=&\bar g(\bar A(\bar e_i,\bar e_j),\bar e_k)
=\bar g(\bar\nabla_{\bar e_j}\bar e_i-\hat{\bar\nabla}_{\bar e_j}\bar e_i,\bar e_k)\nnm\\
=&g(\vfi^{-1}_*(\bar\nabla_{\bar e_j}\bar e_i)-\vfi^{-1}_*(\hat{\bar\nabla}_{\bar e_j}\bar e_i),e_k) =g((\nabla_{e_j}e_i-\hat{\nabla}_{e_j}e_i),e_k) = A(e_i,e_j,e_k).
\end{align}
Consequently
\begin{align}
\bar A=&\sum \bar A(\bar e_i,\bar e_j,\bar e_k)\bar\omega^i\bar\omega^j\bar\omega^k =\sum A(e_i,e_j,e_k)(\vfi^{-1})^*\omega^i(\vfi^{-1})^*\omega^j(\vfi^{-1})^*\omega^k\nnm\\ =&(\vfi^{-1})^*\left(\sum A(e_i,e_j,e_k)\omega^i\omega^j\omega^k\right)=(\vfi^{-1})^*A,
\end{align}
or equivalently, $\vfi^*\bar A=A$. We are done.

Given $c\in\bbr$ and a Riemannian manifold $(M^d,g)$, denote by $\ol{\mathcal S}_{(M^d,g)}(c)$ the set of
all $TM^d$-valued symmetric bilinear forms $A\in \Gamma(\bigodot^2(T^*M^d)\bigotimes (TM^d))$, satisfying the following conditions:

(1) Under the metric $g$, the corresponding $3$-form $A\in \Gamma(\bigodot^2(T^*M^d)\bigotimes (T^*M^d))$ is totally symmetric, that is, $A\in \Gamma(\bigodot^3(T^*M^d))$;

(2) Affine Gauss equation, that is, for any $X,Y,Z\in {\mathfrak X}(M^d)$
\be\label{pre gaus_af sph1}
R(X,Y)Z=c(g(Y,Z)X-g(X,Z)Y)-[A(X),A(Y)](Z).
\ee

Moreover, by adding to $\ol{\mathcal S}_{(M^d,g)}(c)$ the following so called apolarity condition

(3) $\tr_g(A)=0$,

\noindent
we define
$$
{\mathcal S}_{(M^d,g)}(c)=\{A\in\ol{\mathcal S}_{(M^d,g)}(c),\ \tr_g (A)\equiv 0\}.
$$

From Theorem \ref{affine existence} and Theorem \ref{affine uniqueness}, we have
\begin{cor}\label{cor2.1}
For each $A\in {\mathcal S}_{(M^d,g)}(c)$, there uniquely exists one affine hypersphere $x:M^d\to\bbr^{d+1}$ with affine metric $g$, Fubini-Pick form $A$ and affine mean curvature $c$.
\end{cor}

Motivated by Theorem \ref{affine uniqueness}, we introduce the following modified equiaffine equivalence relation between nondegenerate hypersurfaces:

{\dfn Let $x:M^n\to \bbr^{n+1}$ be a nondegenerate hypersurface with the affine metric $g$. A hypersurface $\bar x:M^n\to \bbr^{n+1}$ is called affine equivalent to $x$ if there exists a unimodular transformation $\Phi:\bbr^{n+1}\to \bbr^{n+1}$ and an isometry $\vfi$ of $(M^n,g)$ such that $\bar x=\Phi\circ x\circ\vfi^{-1}$}.

To end this section, we would like to recall the following concept:

{\dfn (\cite{lix13}) A nondegenerate hypersurface $x:M^n\to \bbr^{n+1}$ is called affine symmetric (resp. locally affine symmetric) if

(1) the pseudo-Riemannian manifold $(M^n,g)$ is symmetric (resp. locally symmetric) and therefore $(M^n,g)$ can be written (resp. locally written) as $G/K$ for some connected Lie group $G$ of isometries with $K$ one of its closed subgroups;

(2) the Fubini-Pick form $A$ is invariant under the action of $G$.}

\subsection{The multiple Calabi product of hyperbolic affine hyperspheres}

For later use we make a brief review of the Calabi composition of multiple factors of hyperbolic affine hypersurfaces.
Detailed proofs of the formulas in this subsection has been given in the preprint \cite{lix11}.
\newcommand{\stx}[2]{\strl{(#1)}{#2}}
\newcommand{\spec}[1]{\prod_{#1=1}^K\fr{c_{#1}^{n_{#1}+1}H_{(#1)}^{\fr1{n_{#1}+2}}} {(n_{#1}+1)(-\!\!\stx{#1}{L}_1)}}
\newcommand{\la}{\stx{a}{L}\!\!_1{}}\newcommand{\lb}{\stx{b}{L}\!\!_1{}}
\newcommand{\lc}{\stx{c}{L}\!\!_1{}}\newcommand{\lalp}{\stx{\alpha}{L}\!\!_1{}}
\newcommand{\ha}{\!\stx{a}{h}{}\!\!} %\newcommand{\hb}{\stx{b}{h}\!\!{}}
\newcommand{\Ha}{H_{(a)}}\newcommand{\Hb}{H_{(b)}}\newcommand{\Hc}{H_{(c)}}
\newcommand{\ga}{\!\!\stx{a}{g}{}\!\!\!}\newcommand{\Ga}{\stx{a}{G}\!\!{}}
\newcommand{\gb}{\!\!\stx{b}{g}{}\!\!\!}\newcommand{\Gb}{\stx{b}{G}\!\!{}}
\newcommand{\galp}{\!\!\stx{\alpha}{g}{}\!\!\!}
\newcommand{\xai}{x_{a,i_a}} \newcommand{\xaij}{x_{a,i_aj_a}}
\newcommand{\HH}{f_K\prod_a\fr{c_a^{(n_a+1)(f_K-1)}\Ha^{\fr{f_K+1}{n_a+2}}}
{(n_a+1)^{f_K-n_a}(-\!\!\la)^{f_K-n_a-1}}}
\newcommand{\h}{f^{-\fr1{n+2}}_K \prod_a\fr {(n_a+1)^{\fr{f_K-n_a}{f_K+1}}(-\!\!\la)^{\fr{f_K-n_a-1}{f_K+1}}}
{\left(c_a^{n_a+1}\right)^{\fr{f_K-1}{f_K+1}}\Ha^{\fr1{n_a+2}}}}
\newcommand{\oma}{\stx{a}{\omega}{}\!\!}
\newcommand{\omb}{\stx{b}{\omega}{}\!\!}
\newcommand{\tdca}{(n_a+1)\big(-\!\!\la\big)\prod_b\fr{c_b^{n_b+1}} {(n_b+1)\big(-\!\!\lb\big)}}
\newcommand{\cha}{\prod_a\fr{c_a^{n_a+1}\Ha^{\fr1{n_a+2}}} {(n_a+1)\big(-\!\!\la\big)}}
\newcommand{\chb}{\prod_b\fr{c_b^{n_b+1}\Hb^{\fr1{n_b+2}}} {(n_b+1)\big(-\!\!\lb\big)}}
\newcommand{\chc}{\prod_c\fr{c_c^{n_c+1}\Hc^{\fr1{n_c+2}}} {(n_c+1)\big(-\!\!\lc\big)}}
\newcommand{\ca}{\prod_a\fr{c_a^{n_a+1}}{(n_a+1)\big(-\!\!\la)}}
\newcommand{\cb}{\prod_b\fr{c_b^{n_b+1}}{(n_b+1)\big(-\!\!\lb)}}
\newcommand{\cc}{\prod_c\fr{c_c^{n_c+1}}{(n_c+1)\big(-\!\!\lc)}}
\newcommand{\Aa}{\stx{a}{A}{}\!\!}
\newcommand{\Aalp}{\stx{\alpha}{A}{}\!\!}
\newcommand{\olomea}{\stx{a}{\ol\omega}{}\!\!}
\newcommand{\olomeb}{\stx{b}{\ol\omega}{}\!\!}
\newcommand{\olgma}{\stx{a}{\ol\Gamma}{}\!\!\!}
\newcommand{\olgmb}{\stx{b}{\ol\Gamma}{}\!\!\!}

Now let $r,s$ be two nonnegative integers with $K:=r+s\geq 2$ and $x_\alpha:M^{n_\alpha}_\alpha\to\bbr^{n_\alpha+1}$, $1\leq \alpha\leq s$, be hyperbolic affine hyperspheres of dimension $n_\alpha>0$ with affine mean curvatures $\stx{\alpha}{L}\!\!_1$ and with the origin their common affine center. For convenience we make the following convention:
$$1\leq a,b,c\cdots\leq K,\quad 1\leq\lambda,\mu,\nu\leq K-1,\quad
1\leq\alpha,\beta,\gamma\leq s,\quad \td\alpha=\alpha+r,\ \td\beta=\beta+r,\ \td\gamma=\gamma+r.
$$
Furthermore, for each $\alpha=1,\cdots,s$, set $\td i_{\alpha}=i_\alpha+K-1+\sum_{\beta<\alpha}n_\beta$ with $1\leq i_\alpha\leq n_\alpha$.

Define
$$
f_a:=\begin{cases} a,&1\leq a\leq r;\\ \sum_{\beta\leq \alpha}n_\beta+\td{\alpha},&r+1\leq a=\td\alpha\leq r+s,
\end{cases}
$$
and
$$e_a:=\exp\left(-\fr{t_{a-1}}{n_{a}+1}+\fr{t_{a}}{f_{a}}+\fr{t_{a+1}}{f_{a+1}} +\cdots+\fr{t_{K-1}}{f_{K-1}}\right),\quad 1\leq a\leq K=r+s$$
In particular,
$$
e_1=\exp\left(\fr{t_1}{f_1}+\fr{t_2}{f_2} +\cdots+\fr{t_{K-1}}{f_{K-1}}\right),\quad
e_K=\exp\left(-\fr{t_{K-1}}{n_K+1}\right).
$$

Put $n=\sum_\alpha n_\alpha+K-1$ and $M^n=R^{K-1}\times M^{n_1}_1\times\cdots\times M^{n_s}_s$. For any $K$ positive numbers $c_1,\cdots,c_K$, define a smooth map $x:M^n\to\bbr^{n+1}$ by
\begin{align}
x(t^1,&\cdots,t^{K-1},p_1,\cdots,p_s):=(c_1e_1,\cdots, c_re_r,c_{r+1}e_{r+1}x_1(p_1),\cdots,c_Ke_Kx_s(p_s)),\nnm\\&\hs{1cm}\forall (t^1,\cdots,t^{K-1},p_1,\cdots,p_s)\in M^n.\label{mulpro2}
\end{align}

{\prop\label{general sense} (\cite{lix11}) The map $x:M^n\to\bbr^{n+1}$ defined above is a new hyperbolic affine hypersphere with the affine mean curvature
\be\label{newl1c}
L_1=-\fr1{(n+1)C},\quad C:=\left(\fr1{n+1}\prod_{a=1}^r c_a^2\cdot\prod_{\alpha=1}^s\fr{c_{r+\alpha}^{2(n_\alpha+1)}} {(n_\alpha+1)^{n_\alpha+1}(-\!\!\stx{\alpha}{L}_1)^{n_\alpha+2}}\right)^{\fr1{n+2}},
\ee
Moreover, for given positive numbers $c_1,\cdots,c_K$, there exits some $c>0$ and $c'>0$ such that
the following three hyperbolic affine hyperspheres
\bea &x:=(c_1e_1,\cdots, c_re_r,c_{r+1}e_{r+1}x_1,\cdots,c_Ke_sx_s),\nnm\\
&\bar x:=c(e_1,\cdots, e_r,e_{r+1}x_1,\cdots,e_sx_s),\nnm\\
&\td x:=(e_1,\cdots, e_r,e_{r+1}x_1,\cdots,c'e_sx_s)\nnm
\eea
are equiaffine equivalent to each other.}

{\dfn\label{df2}(\cite{lix11}) \rm The hyperbolic affine hypersphere $x$ is called the Calabi composition of $r$ points and $s$ hyperbolic affine hyperspheres.}

Denote by $\{v^{i_\alpha}_\alpha; \ i_\alpha=1,\cdots,n_\alpha\}$ the local coordinate system of $M_\alpha$, $\alpha=1,\cdots,s$. Then we have

{\prop\label{corr0}(\cite{lix11}) The affine metric $g$, the affine mean curvature $L_1$ and the possibly nonzero components of the Fubini-Pick form $A$ of the Calabi composition $x:M^n\to \bbr^{n+1}$ of $r$ points and $s$ hyperbolic affine hyperspheres $x_\alpha:M_\alpha\to\bbr^{n_\alpha+1}$, $\alpha=1,\cdots,s$, are given as follows:
\begin{align}
&g_{\lambda\mu}=\begin{cases} \displaystyle\fr{\lambda+1}{\lambda}C\delta_{\lambda\mu},&\!1\leq\lambda\leq r-1;\\
\displaystyle\fr{n_1+r+1}{r(n_1+1)}C\delta_{r\mu},&\!\lambda=r;\\
\displaystyle\fr{\sum_{\beta\leq\alpha+1}n_\beta+\td{\alpha}+1}  {(n_{\alpha+1}+1)(\sum_{\beta\leq \alpha}n_\beta+\td{\alpha})}C\delta_{\lambda\mu}, &\!\!r+1\leq\lambda=\td\alpha\leq r+s-1.\end{cases}\label{g-lam mu}
\\
&g_{\td i_{\alpha}\td j_{\beta}}=(n_\alpha+1) (-\!\!\lalp)C\galp_{i_{\alpha}j_{\alpha}}\delta_{\alpha\beta},
\quad
g_{\lambda\td i_{\td\alpha}}=0.\label{g-gen}
\\
&A_{\lambda\lambda\lambda} =\begin{cases}\displaystyle \fr{1-\lambda^2}{\lambda^2}C,&1\leq\lambda\leq r-1,\\
\displaystyle\left(\fr1{r^2}-\fr1{(n_1+1)^2}\right)C,&\lambda=r,\\
\displaystyle\fr{(\sum_{\beta\leq\alpha+1}n_\beta+\td{\alpha}+1)C} {(n_{\alpha+1}+1)(\sum_{\beta\leq \alpha}n_\beta+\td{\alpha})}\left(\fr1{\sum_{\beta\leq \alpha}n_\beta+\td{\alpha}}-\fr1{n_{\alpha+1}+1}\right), &r+1\leq\lambda=\td\alpha\leq r+s-1.
\end{cases}
\\
&A_{\lambda\lambda\mu} =\begin{cases}\displaystyle \fr{\lambda+1}{\lambda\mu}C,&1\leq\lambda<\mu\leq r,\\
\displaystyle\fr{(\lambda+1)C}{\lambda(\sum_{\beta\leq \alpha}n_\beta+\td{\alpha})},&1\leq\lambda\leq r-1, \mu=\td\alpha,\\
\displaystyle\fr{(n_1+r+1)C}{r(\sum_{\beta\leq \alpha}n_\beta+\td{\alpha})},&\lambda=r,\ \mu=\td\alpha,\\
\displaystyle\fr{(\sum_{\gamma\leq \alpha+1}n_\gamma+\td{\alpha}+1)C} {(n_{\alpha+1}+1)(\sum_{\gamma\leq \alpha}n_\gamma+\td{\alpha})(\sum_{\gamma\leq \beta}n_\gamma+\td{\beta})},&r+1\leq\lambda=\td\alpha<\mu=\td\beta\leq r+s-1.
\end{cases}
\\
&A_{\td i_{\alpha}\td j_{\alpha}\,{\td\alpha}-1} =-\fr1{n_\alpha+1}g_{\td i_{\alpha}\td j_{\alpha}}=-\!\lalp C\galp_{i_\alpha j_\alpha},\\
&A_{\td i_{\alpha}\td j_{\alpha}\td\beta}=\fr1{\sum_{\gamma\leq \beta}n_\gamma+\td{\beta}}g_{\td i_{\alpha}\td j_{\alpha}}=\fr{(n_\alpha+1)\big(-\!\!\lalp\big)C}{\sum_{\gamma\leq \beta}n_\gamma+\td\beta}\ \galp_{i_\alpha j_\alpha},\quad \beta\geq \alpha,\\
&A_{\td i_{\alpha}\td j_{\alpha}\td k_{\alpha}}=(n_\alpha+1)\big(-\!\!\lalp\big)C\Aalp_{i_\alpha j_\alpha k_\alpha},
\end{align}
where $\stx{\alpha}{L}_1$, $\stx{\alpha}{g}$ and $\stx{\alpha}A$ are the affine mean curvature, the affine metric and the Fubini-Pick form of $x_\alpha$, $\alpha=1,\cdots,s$.}

From Proposition \ref{corr0}, the following corollary is easily derived (cf. \cite{lix11}):

{\cor The Calabi composition $x:M^n\to \bbr^{n+1}$ of $r$ points and $s$ hyperbolic affine hyperspheres $x_\alpha:M^{n_\alpha}\to\bbr^{n_\alpha+1}$, $1\leq \alpha\leq s$, is affine symmetric if and only if
each positive dimensional factor $x_\alpha$ is symmetric.}

The next example will be used later:

\expl\label{expl}\rm Given a positive number $C_0$, let $x_0:\bbr^{n_0}\to \bbr^{n_0+1}$ be the well known flat hyperbolic affine hypersphere of dimension $n_0$ which is defined by
$$
x^1\cdots x^{n_0} x^{n_0+1}=C_0,\quad x^1>0,\cdots,x^{n_0+1}>0.
$$
Then it is not hard to see that $x_0$ is the Calabi composition of $n_0+1$ points. In fact, we can write for example
$$
x_0=(e_1,\cdots,e_{n_0},C_0e_{n_0+1}).
$$
Then by Corollary \ref{corr0} the affine metric $g_0$, the affine mean curvature $\stx{0}{L}_1$ and the Fubini-Pick form $\stx{0}{A}$ of $x_0$ are respectively given by (cf. \cite{li-sim-zhao93})
\begin{align}
\stx{0}{g}_{\lambda\mu}=& \fr{\lambda+1}{\lambda}\left(\fr{C_0^2}{n_0+1}\right)^{\fr1{n_0+2}} \delta_{\lambda\mu},\label{gofexpl}\\
\stx{0}{L}_1=&-\fr1{(n_0+1)C} =-(n_0+1)^{-\fr{n_0+1}{n_0+2}}C_0^{-\fr2{n_0+2}},\label{l1ofexpl}\\
\stx{0}{A}_{\lambda\mu\nu}=&\begin{cases}-\fr{\lambda^2-1}{\lambda^2} \left(\fr{C_0^2}{n_0+1}\right)^{\fr1{n_0+2}}, &{\rm if\ }\lambda=\mu=\nu;\\
\fr{\lambda+1}{\lambda\nu}\left(\fr{C_0^2}{n_0+1}\right)^{\fr1{n_0+2}}, &{\rm if\ }\lambda=\mu<\nu;\\
0, &{\rm otherwise.}
\end{cases}\label{aofexpl}
\end{align}
Thus the Pick invariant of $x_0$ is
\be\label{jofexpl}
\stx{0}{J}=\fr1{n_0(n_0-1)}\stx{0}{g}{}\!\!^{\lambda_1\lambda_2} \stx{0}{g}{}\!\!^{\mu_1\mu_2} \stx{0}{g}{}\!\!^{\nu_1\nu_2} \stx{0}{A}_{\lambda_1\mu_1\nu_1}\stx{0}{A}_{\lambda_2\mu_2\nu_2} =(n_0+1)^{-\fr{n_0+1}{n_0+2}}C_0^{-\fr2{n_0+2}}=-\!\!\stx{0}{L}_1.
\ee

By restrictions, $g$ defines a flat metric $g_0$ on $\bbr^{K-1}$ with matrix $(g_{\lambda\mu})$ and, for each $\alpha$, a metric $g_\alpha$ on $M_\alpha$ with matrix $\big(g^\alpha_{i_\alpha j_\alpha}\big)=\big(g_{\td i_{\alpha}\td j_{\alpha}}\big)$ and inverse matrix $\big(g^{i_\alpha j_\alpha}_\alpha\big)$, which is conformal to the original metric $\stx{\alpha}{g}$, or more precisely, \be\label{g_alpha}g_\alpha=(n_\alpha+1) \big(-\!\!\stx{\alpha}{L}_1)C\stx{\alpha}{g}.\ee

Write $M_0=\bbr^{K-1}$. Then, with respect to the affine metric $g$ on $M^n$, the Fubini-Pick form $A$ can be identified with a $TM^n$-valued symmetric $2$-form $ A:TM^n\times TM^n\to TM^n$. For each ordered triple $\alpha,\beta,\gamma\in\{0,1,\cdots,s\}$, $ A$ defines one $TM_\gamma$-valued bilinear map $ A^{\gamma}_{\alpha\beta}:TM_{\alpha}\times TM_{\beta}\to TM_\gamma$, which is the $TM_\gamma$-component of $ A_{\alpha\beta}$, the restriction of $ A$ to $TM_{\alpha}\times TM_{\beta}$. For $\alpha=1,\cdots,s$, define
\be\label{Halpha}
H_\alpha=\fr1{n_{\alpha}}\tr_{g_\alpha} A^0_{\alpha\alpha}\equiv \fr1{n_{\alpha}}g^{i_\alpha j_\alpha}_\alpha A^0_{\alpha\alpha}\left(\pp{}{v^{i_{\alpha}}_{\alpha}}, \pp{}{v^{j_{\alpha}}_{\alpha}}\right),
\ee
where the metric $g_\alpha$ is given by \eqref{g_alpha}. Then we have (\cite{lix11})

\begin{prop}\label{properties} Let $x:M^n\to \bbr^{n+1}$ be the Calabi composition of $r$ points and $s$ hyperbolic affine hyperspheres and $g$ the affine metric of $x$. Then

$(1)$ The Riemannian manifold $M^n\equiv (M^n,g)$ is reducible, that is
\be\label{derham}
(M^n,g)=\bbr^q\times(M_1,g_1)\times\cdots\times(M_s,g_s),\quad q+s\geq 2;\ee

$(2)$ There must be a positive dimensional Euclidean factor $\bbr^q$ in the de Rham decomposition \eqref{derham} of $M^n$, that is, $q>0$;

$(3)$ $q\geq s-1$ with the equality holding if and only if $r=0$;

$(4)$ $ A^\gamma_{\alpha\beta}\equiv 0$ if
$(\alpha,\beta,\gamma)$ is not one of the following triples: $(0,0,0)$, $(\alpha,\alpha,0)$, $(\alpha,0,\alpha)$, $(0,\alpha,\alpha)$ or $(\alpha,\alpha,\alpha)$.

$(5)$ For any $p=(p_0,p_1,\cdots,p_s)\in M^n$ and each $\alpha=1,\cdots,s$, it holds that
\begin{align}\label{added f0}
A^0_{\alpha\alpha}(R^{M_\alpha}(X_\alpha,Y_\alpha)&Z_\alpha,W_\alpha) +A^0_{\alpha\alpha}(Z_\alpha,R^{M_\alpha} (X_\alpha,Y_\alpha)W_\alpha)=0,\nnm\\
&\forall X_\alpha,Y_\alpha,Z_\alpha,W_\alpha\in T_{p_\alpha}M_\alpha
\end{align}
equivalent to that the holonomy algebra ${\mathfrak h}_\alpha$ of $(M_\alpha,g_\alpha)$ acts on $A^0_{\alpha\alpha}$ trivially, that is
${\mathfrak h}_\alpha\cdot A^0_{\alpha\alpha}=0$.

$(6)$ The vector-valued functions $H_\alpha$, $\alpha=1,\cdots,s$, defined by \eqref{Halpha} satisfy the following equalities:
\begin{align}
&H_\alpha=-\fr{f_{\td\alpha-1}}{f_{\td\alpha}C}\pp{}{t^{\td\alpha-1}} +\fr1C\sum_{s-1\geq\beta\geq\alpha}\fr{n_{\beta+1}+1} {f_{\td\beta+1}}\pp{}{t^{\td\beta}},\label{Halpha1}\\
&g(H_\alpha,H_\alpha)=C^{-1}\left(\fr1{n_{\alpha}+1}-\fr1{f_K}\right) =\fr{n-n_\alpha}{n_\alpha+1}(-L_1),\\
&g(H_\alpha,H_\beta)=L_1\quad{\rm for\ }\alpha\neq\beta;
\end{align}

$(7)$ $A^{\alpha}_{\alpha\alpha}$ is identical to the $TM_\alpha$-valued symmetric bilinear form defined by the Fubini-Pick form $\stx{\alpha}{A}$ of $x_\alpha$.
\end{prop}

In the next section we shall prove as the main result that a locally strongly convex affine hypersurface $x:M^n\to \bbr^{n+1}$ is locally the Calabi composition of some points and hyperbolic affine hyperspheres if and only if the above conditions (1), (4) and (5) hold.

\section{A characterization of the Calabi composition}

We are mainly to establish a necessary and sufficient condition for a locally strongly convex hyperbolic hypersphere locally to be the Calabi composition of several hyperbolic affine hyperspheres, possibly including point factors. It turns out that
this special characterization theorem is needed in some related important classifications.

\begin{thm}[The Main Theorem]\label{the main} A locally strongly convex hyperbolic hypersphere $x:M^n\to\bbr^{n+1}$, with the affine metric $g$ and the Fubini-Pick form $ A\in\Gamma(T*M^n\bigodot T^*M^n\bigotimes TM^n)$, is locally affine equivalent to the Calabi composition of some hyperbolic affine hyperspheres, possibly including point factors, if and only if the following three conditions hold:

(1) The Riemannian manifold $(M^n,g)$ is reducible, that is
\be\label{derham1}
(M^n,g)=\bbr^q\times(M_1,g_1)\times\cdots\times(M_s,g_s),\quad q+s\geq 2,\ee
where $(M_\alpha,g_\alpha)$, $\alpha=1,\cdots,s$, are irreducible Riemannian manifolds;

(2) Denote $M_0=\bbr^q$ and by $ A^{\gamma}_{\alpha\beta}$, $0\leq\alpha,\beta,\gamma\leq s$, the $TM_\gamma$-component of the restriction of $ A$ to $TM_{\alpha}\times TM_{\beta}$. Then
$ A^\gamma_{\alpha\beta}\equiv 0$ if the oriented triple
$(\alpha,\beta,\gamma)$ is not one of the followings:
$$
(0,0,0),\quad (\alpha,\alpha,0),\quad (\alpha,0,\alpha),\quad (0,\alpha,\alpha),\quad (\alpha,\alpha,\alpha).
$$

(3) For each $\alpha=1,\cdots,s$, the holonomy algebra ${\mathfrak h}_\alpha$ of $(M_\alpha,g_\alpha)$ acts on $A^0_{\alpha\alpha}$ trivially:
${\mathfrak h}_\alpha\cdot A^0_{\alpha\alpha}=0$, or equivalently, for any $p=(p_0,p_1,\cdots,p_s)\in M^n$ and for all $X_\alpha,Y_\alpha,Z_\alpha,W_\alpha\in T_{p_\alpha}M_\alpha$, it holds that
\be\label{halpha a0alpha}
A^0_{\alpha\alpha}(R^{M_\alpha}(X_\alpha,Y_\alpha)Z_\alpha,W_\alpha)+A^0_{\alpha\alpha}(Z_\alpha,R^{M_\alpha}(X_\alpha,Y_\alpha)W_\alpha)=0.
\ee
\end{thm}

\proof
Since our consideration is local here, we suppose that $M^n$ is connected and simply connected. The necessary part of the theorem is clearly from Proposition \ref{properties}. To prove the sufficient part, we first note that, the affine mean curvature $L_1$ of a hyperbolic affine hypersphere is a negative constant.

{\bf Claim}:
$M^n$ must have an Euclidean factor $\bbr^{q}$, $q>0$, in its de Rham decomposition \eqref{derham1}.

In fact it suffices to show that $M^n$ would be irreducible if $q=0$. To do this we suppose that $s>1$. Then for non-vanishing $X\in TM_1$ and $Y\in TM_2$, we have by the affine Gauss equation \eqref{gaus_af sph} that
\begin{align}
0=& R(X,Y)Y=L_1(g(Y,Y)X-g(X,Y)Y)\nnm\\
&\hs{2.5cm}-[A(X),A(Y)](Y)=L_1g(Y,Y)X,\label{temp}
\end{align}
where we have used the fact that $A(Y)(Y)=A(Y,Y)\in TM_2$ and $A(X)(Y)=A(X,Y)=0$ due to the condition (2) in the theorem. Clearly \eqref{temp} is not possible since $L_1<0$ and $Y\neq 0$. This contradiction proves the claim.

By the assumption, $ A$ can be decomposed as
\be\label{decomposition of A}
 A=\sum_{\alpha=0}^s A_{\alpha\alpha}^\alpha +\sum_{\alpha=1}^s A_{\alpha\alpha}^0 +\sum_{\alpha=1}^s A_{\alpha 0}^\alpha+\sum_{\alpha=1}^s A_{0\alpha}^\alpha.
\ee

Same as in the last section, for $\alpha=1,\cdots,s$, we define $H_\alpha=\fr1{n_\alpha}\tr_{ g_\alpha}( A^0_{\alpha\alpha})$ and denote $ \bar c_\alpha=| H_\alpha|$. Then we have
\begin{lem}\label{lem added}
$A^0_{00}\in \ol{\mathcal S}_{\bbr^q}(L_1)$ and the following identities hold:
\begin{align}
&A^0_{\alpha\alpha}(X_\alpha,Y_\alpha)=g(X_\alpha,Y_\alpha)H_\alpha,\quad g_0(H_\alpha,H_\beta)=L_1,\nnm\\
&\hs{1.5cm}{\rm if\ }1\leq\alpha\neq\beta\leq s,\label{alp alp 0}\\
& A^\alpha_{\alpha 0}(X_\alpha,Z_0) = A^\alpha_{0\alpha}(Z_0,X_\alpha)=g_0(Z_0,H_\alpha)X_\alpha.\label{alp 0 alp}\\
& A^0_{00}(Z_0,H_\alpha) =g_0(Z_0,H_\alpha)H_\alpha+L_1Z_0,\label{000}
\end{align}
where $Z_0\in T_{p_0}M_0$, $X_\alpha,Y_\alpha\in T_{p_\alpha}M_\alpha$, and $g_0$ is the flat metric on $\bbr^q$.
\end{lem}

{\it Proof of Lemma \ref{lem added}}: First
note that the first conclusion is direct from the fact that $A\in {\mathcal S}_{(M^n,g)}(L_1)$ together with the decomposition \eqref{decomposition of A}.

Let $p=(p_0,p_1,\cdots,p_s)\in M^n$ be an arbitrary point with $p_\alpha\in M_\alpha$, $\alpha=0,1,\cdots,s$. Let ${\mathfrak h}$ be the holonomy algebra of $(M^n,g)$, and ${\mathfrak h}_1$, $\cdots$, ${\mathfrak h}_s$ the holonomy algebras of $(M_1,g_1)$, $\cdots$, $(M_s,g_s)$, respectively. Then by \eqref{derham1} we have
\be\label{holo decom}
{\mathfrak h}={\mathfrak h}_1+\cdots+{\mathfrak h}_s,
\ee
Since Riemannian manifold $(M_\alpha,g_\alpha)$ ($1\leq\alpha\leq s$) is irreducible, ${\mathfrak h}_\alpha$ acts irreducible on $T_{p_\alpha}M_\alpha$ for every point $p_\alpha\in M_\alpha$ and, at the same time, acts trivially on any other $T_{p_\beta}M_\beta$, $\beta\neq\alpha$.

For each $\alpha=1,\cdots,s$ and for all $T\in{\mathfrak h}_\alpha$, $X_\alpha,Y_\alpha\in T_{p_\alpha}M_\alpha$,
\eqref{halpha a0alpha} gives that
\be\label{added f1}
A^0_{\alpha\alpha}(TX_\alpha,Y_\alpha) +A^0_{\alpha\alpha}(X_\alpha,TY_\alpha)=0.
\ee
By using the irreducibility of ${\mathfrak h}_\alpha$ on $T_{p_\alpha}M_\alpha$, we get from \eqref{added f1}
\be
\lagl A^0_{\alpha\alpha}(X_\alpha,Y_\alpha),e_a\ragl=c^a_\alpha g_\alpha(X_\alpha,Y_\alpha)
\ee
for some constants $c^a_\alpha\in\bbr$, where $\{e_a\}$ is an orthonormal basis for $T_{p_0}\bbr^{q}\equiv\bbr^q$.
Therefore
$$
A^0_{\alpha\alpha}(X_\alpha,Y_\alpha) =g_\alpha(X_\alpha,Y_\alpha)\sum_ac^a_\alpha e_a.
$$

On the other hand, it is seen that $H_\alpha=\sum_ac^a_\alpha e_a$, proving the first equality in \eqref{alp alp 0}. This with the symmetry of $A$ gives \eqref{alp 0 alp}.

Now put $Z_0\in T_{p_0}M_0$, $X_\alpha,Y_\alpha\in T_{p_\alpha}M_\alpha$. Then the affine Gauss equation \eqref{gaus_af sph} of $x$ tells that
\be\label{added f2}
0=L_1g_\alpha(X_\alpha,Y_\alpha)Z_0-[A(Z_0),A(X_\alpha)](Y_\alpha).
\ee
But by the first equality of \eqref{alp alp 0} and \eqref{alp 0 alp} we find
\begin{align*}
&[A(Z_0),A(X_\alpha)](Y_\alpha)=A(Z_0,A(X_\alpha,Y_\alpha)) -A(X_\alpha,A(Z_0,Y_\alpha))\\ =&A(Z_0,A^0_{\alpha\alpha}(X_\alpha,Y_\alpha) +A^\alpha_{\alpha\alpha}(X_\alpha,Y_\alpha)) -A(X_\alpha,A^\alpha_{0\alpha}(Z_0,Y_\alpha))\\
=&A(Z_0,g_\alpha(X_\alpha,Y_\alpha)H_\alpha) +A(Z_0,A^\alpha_{\alpha\alpha}(X_\alpha,Y_\alpha)) -A(X_\alpha,g(Z_0,H_\alpha)Y_\alpha))\\
=&g_\alpha(X_\alpha,Y_\alpha)A^0_{00}(Z_0,H_\alpha) +A^\alpha_{0\alpha}(Z_0,A^\alpha_{\alpha\alpha}(X_\alpha,Y_\alpha)) \\ &\ -g(Z_0,H_\alpha)(A^0_{\alpha\alpha}(X_\alpha,Y_\alpha) +A^\alpha_{\alpha\alpha}(X_\alpha,Y_\alpha))\\
=&g_\alpha(X_\alpha,Y_\alpha)A^0_{00}(Z_0,H_\alpha)+g(Z_0,H_\alpha) A^\alpha_{\alpha\alpha}(X_\alpha,Y_\alpha)\\
&\ -g(Z_0,H_\alpha)g_\alpha(X_\alpha,Y_\alpha)H_\alpha -g(Z_0,H_\alpha)A^\alpha_{\alpha\alpha}(X_\alpha,Y_\alpha)\hs{2cm}\mb{}\\
=&g_\alpha(X_\alpha,Y_\alpha)(A^0_{00}(Z_0,H_\alpha) -g(Z_0,H_\alpha)H_\alpha).
\end{align*}
Putting this equality into \eqref{added f2} we obtain \eqref{000}.

Similarly, for $1\leq\alpha\neq\beta\leq s$, let $X_\alpha\in T_{p_\alpha}M_\alpha$ and $Y_\beta,Z_\beta\in T_{p_\beta}M_\beta$. Then
\be\label{added f3}
0=L_1g_{\beta}(Y_\beta,Z_\beta)X_\alpha-[A(X_\alpha),A(Y_\beta)](Z_\beta).
\ee

But
\begin{align*}
[A(X_\alpha),A(Y_\beta)](Z_\beta)=&A(X_\alpha,A(Y_\beta,Z_\beta)) -A(Y_\beta,A(X_\alpha,Z_\beta))\\ =&A(X_\alpha,A^0_{\beta\beta}(Y_\beta,Z_\beta) +A^\beta_{\beta\beta}(Y_\beta,Z_\beta))\\
=&A^\alpha_{\alpha0}(X_\alpha,g_\beta(Y_\beta,Z_\beta)H_\beta) =g_\beta(Y_\beta,Z_\beta)g(H_\alpha,H_\beta)X_\alpha.
\end{align*}
Comparing this with \eqref{added f3} gives the second formula in \eqref{alp alp 0}.
\endproof

For the given point $p=(p_0,p_1,\cdots,p_s)\in M^n$, denote by ${\mathcal H}_p$ the subspace of $T_{p_0}M_0\equiv \bbr^q$ generated by
$H_1(p_1),\cdots H_s(p_s)$ and by
${\mathcal H}^\bot_p$ the orthogonal complement of ${\mathcal H}_p$ in $T_{p_0}M_0$. Thus $T_{p_0}M_0={\mathcal H}_p\oplus {\mathcal H}^\bot_p$. Then $A^0_{00}(p)$ can be decomposed into the sum of its ${\mathcal H}_p$-component $ A_0^{\mathcal H_p}$ and its ${\mathcal H^\bot_p}$-component $ A_0^{\mathcal H^\bot_p}$, that is,
$ A^0_{00}(p)= A_0^{\mathcal H_p}+ A_0^{\mathcal H^\bot_p}$.

\begin{lem}\label{rh geq s-1} Let $q$, $s$ be as above. Then
$q\geq \dim {\mathcal H_p}\geq s-1$. Furthermore, $q\geq s$ if and only if $\dim {\mathcal H_p}=s$.\end{lem}

{\it Proof of Lemma \ref{rh geq s-1}}: For simplicity we omit the point $p$ in the symbols.
To prove the first part of the lemma, it suffices to show that the set of the $s$ nonzero vectors $H_1,\cdots,H_s$ in $TM_0\equiv\bbr^{q}$ has a rank not less than $s-1$. This is equivalent to show that the $s$-th order matrix
$$
\lmx g(H_1,H_1)&g(H_1,H_2)&\cdots&g(H_1,H_s)\\
g(H_2,H_1)&g(H_2,H_2)&\cdots&g(H_2,H_s)\\
\cdots&\cdots&\cdots&\cdots\\
g(H_s,H_1)&g(H_s,H_2)&\cdots&g(H_s,H_s)
\rmx
=\lmx  \bar c^2_1&L_1&\cdots&L_1\\
L_1& \bar c^2_2&\cdots&L_1\\
\cdots&\cdots&\cdots&\cdots\\
L_1&L_1&\cdots& \bar c^2_s
\rmx
$$
has a rank equal to or larger than $s-1$.

Indeed, by deleting the last line and the second last column, we find a ($s-1$)-minor of the above matrix:
\begin{align}
&\det\lmx  \bar c^2_1&L_1&\cdots&L_1\\
L_1& \bar c^2_2&\cdots&L_1\\
\cdots&\cdots&\cdots&\cdots\\
L_1&L_1&\cdots&L_1
\rmx
=L_1\det\lmx  \bar c^2_1&L_1&\cdots&L_1\\
L_1& \bar c^2_2&\cdots&L_1\\
\cdots&\cdots&\cdots&\cdots\\
1&1&\cdots&1
\rmx\nnm\\
=&L_1\det\lmx  \bar c^2_1-L_1&0&\cdots&0\\
0& \bar c^2_2-L_1&\cdots&0\\
\cdots&\cdots&\cdots&\cdots\\
1&1&\cdots&1
\rmx\nnm\\
=&L_1( \bar c^2_1-L_1)\cdots( \bar c^2_{s-2}-L_1)<0.
\end{align}

Furthermore, if $q\geq s$, and $\dim {\mathcal H}=s-1$, then $r-1:=\dim{\mathcal H}^\bot\geq 1$. Consider the restriction $\bar A_0^{\mathcal H}$ of $  A_0^{\mathcal H}$ to the subspace ${\mathcal H}^\bot\times{\mathcal H}^\bot$. Define $H_0=\fr1{r-1}\tr\bar A_0^{\mathcal H}$ and $ \bar c_0=|H_0|$. Then, for any unit vector $e_0\in {\mathcal H}^\bot$ and each $\alpha=1,\cdots,s$, we have
\begin{align*}
g(\bar A^{\mathcal H}_0(e_0,e_0),H_\alpha)=&g( A(e_0,e_0),H_\alpha) =g( A^0_{00}(e_0,H_\alpha),e_0)\\
=&L_1g(e_0,e_0)=L_1,
\end{align*}
implying that
$$
g(H_0,H_\alpha)=L_1,\quad \alpha=1,\cdots,s.
$$

Then in the same way as in proving that the rank of the matrix $(g(H_\alpha,H_\beta))_{1\leq\alpha,\beta\leq s}$ of order $s$ is no less than $s-1$ we can obtain that the rank of the $(s+1)$-th order matrix $(g(H_\alpha,H_\beta))_{0\leq\alpha,\beta\leq s}$ is no less than $s$. Since $\{H_\alpha;\ 0\leq\alpha\leq s\}\subset {\mathcal H}$, it follows that $\dim {\mathcal H}\geq s$ which contradicts the assumption.\endproof

Since $q$ is fixed, we have
\begin{cor}
$\dim {\mathcal H_p}$ and $\dim {\mathcal H^\bot_p}$ are independent of the point $p$. Thus they form two subbundles ${\mathcal H}$ and ${\mathcal H^\bot}$ of the tangent bundle $TM_0$.
\end{cor}

Define $\bar A_0^{\mathcal H^\bot}=\left. A_0^{\mathcal H^\bot}\right|_{{\mathcal H^\bot}\times {\mathcal H^\bot}}$. By \eqref{000}, for any $X_{\mathcal H},Y_{\mathcal H}\in {\mathcal H}$, $ A^0_{00}(X_{\mathcal H},Y_{\mathcal H})\in {\mathcal H}$, and for any $Y_{\mathcal H^\bot}\in {\mathcal H^\bot}$, $ A^0_{00}(X_{\mathcal H},Y_{\mathcal H^\bot})\in {\mathcal H^\bot}$. Therefore, $ A^0_{00}$ can be decomposed into the following components:
\be\label{decomp tdsigma0}
 A^0_{00}=\bar A_0^{\mathcal H^\bot}+\bar A_0^{\mathcal H}+\left. A_0^{\mathcal H^\bot}\right|_{{\mathcal H^\bot}\times {\mathcal H}}+\left. A_0^{\mathcal H^\bot}\right|_{{\mathcal H}\times {\mathcal H^\bot}}+\left. A_0^{\mathcal H}\right|_{{\mathcal H}\times {\mathcal H}}.
\ee
\begin{lem}\label{barsigma0} Define $r=\rank{\mathcal H}^\bot+1$. Then $\bar A_0^{\mathcal H^\bot}\in {\mathcal S}_{\bbr^{r-1}}(\fr{(n+1)L_1}r)$ if $\rank{\mathcal H}^\bot\geq 1$.\end{lem}

{\it Proof of Lemma \ref{barsigma0}}:
Since $\rank{\mathcal H}^\bot\geq 1$, it holds by Lemma \ref{rh geq s-1} that
$$
q=\rank{\mathcal H}+\rank{\mathcal H}^\bot\geq s-1+1=s.
$$
Making use of Lemma \ref{rh geq s-1} once again we have that $\rank{\mathcal H}=s$ and therefore $\{H_1,\cdots,H_s\}$ is a frame for the vector bundle ${\mathcal H}$. Set $h_{\alpha\beta}=g(H_\alpha,H_\beta)$, $1\leq\alpha,\beta\leq s$, and $(h^{\alpha\beta})=(h_{\alpha\beta})^{-1}$.

We first compute $\bar A^{\mathcal H}_0$. Write
$$
\bar A^{\mathcal H}_0(X,Y)=\sum C^\alpha_{XY}H_\alpha,\quad \forall X,Y\in {\mathcal H^\bot}.
$$

Then we have
\begin{align}
g(\bar A^{\mathcal H}_0(X,Y),H_\alpha) =&\sum C^\beta_{XY}g(H_\beta,H_\alpha)=\sum C^\beta_{XY}h_{\beta\alpha};\nnm\\
g(\bar A^{\mathcal H}_0(X,Y),H_\alpha) =&g( A^0_{00}(X,Y),H_\alpha) =g( A^0_{00}(X,H_\beta),Y)\nnm\\
=&g(g(X,H_\alpha)H_\alpha+L_1X,Y)=L_1g(X,Y),\nnm
\end{align}
implying that
$
\sum C^\beta_{XY}h_{\beta\alpha}=L_1g(X,Y)
$
or equivalently $$C^\alpha_{XY}=L_1g(X,Y)\sum_{\beta} h^{\alpha\beta}.$$
It follows that
\be\label{barsigma0h}
\bar A^{\mathcal H}_0(X,Y)=L_1g(X,Y)\sum_{\alpha,\beta} h^{\alpha\beta}H_\alpha,
\quad\forall X,Y\in {\mathcal H^\bot}.
\ee

Thus we have
\be\label{H0}
H_0=\fr1{r-1}\tr(\bar A^{\mathcal H}_0)=L_1\sum_{\alpha,\beta} h^{\alpha\beta}H_\alpha.
\ee
On the other hand, by using \eqref{decomposition of A} and the fact that $\tr( A)=0$, it is seen that
$$
\tr A^0_{00}+\sum_\alpha\tr A^0_{\alpha\alpha}=0,
$$
which with the decomposition \eqref{decomp tdsigma0} gives
\begin{align}
&\tr(\bar A^{\mathcal H^\bot}_0)=0,\label{tr barsigma0hbot=0}\\
&(r-1)H_0+\sum_{\alpha,\beta}h^{\alpha\beta} A^{\mathcal H}_0(H_\alpha,H_\beta) +\sum_\alpha n_\alpha H_\alpha=0.\label{sum Halpha=0}
\end{align}

But by \eqref{000},
$$
 A^{\mathcal H}_0(H_\alpha,H_\beta)=g(H_\alpha,H_\beta)H_\beta +L_1H_\alpha =h_{\alpha\beta} H_\beta +L_1H_\alpha.
$$
Thus \eqref{sum Halpha=0} can be rewritten as
\be\label{sum Halpha=0-1}
(r-1)H_0+\sum_\alpha(1+L_1\sum_\beta h^{\alpha\beta})H_\alpha +\sum_\alpha n_\alpha H_\alpha=0.
\ee

Comparing \eqref{H0} and \eqref{sum Halpha=0-1} gives that
$$
\sum_\alpha\left(rL_1\sum_\beta h^{\alpha\beta}+(n_\alpha+1)\right)H_\alpha=0,
$$
or equivalently
\be\label{sumbeta h^alphabeta}
rL_1\sum_\beta h^{\alpha\beta}+(n_\alpha+1)=0,\quad \alpha=1,\cdots,s.
\ee

It follows that
\be\label{*}
\sum_\beta h^{\alpha\beta}=-\fr{n_\alpha+1}{rL_1},\quad \forall \alpha
\ee
which with \eqref{barsigma0h} gives
\be\label{barsigma0h-1}
\bar A^{\mathcal H}_0(X,Y)=-g(X,Y)\sum_\alpha \fr{n_\alpha+1}{r}H_\alpha,
\quad\forall X,Y\in {\mathcal H^\bot},
\ee
and thus
\be\label{H0-1}
H_0=\fr1{r-1}\tr(\bar A^{\mathcal H}_0)=-\sum_\alpha \fr{n_\alpha+1}{r}H_\alpha.
\ee

Since we have shown that $\tr(\bar A^{\mathcal H^\bot}_0)=0$ (Equation \eqref{tr barsigma0hbot=0}), to complete the proof of Lemma \ref{barsigma0}, it now suffices to show that
\begin{align}
\fr{(n+1)L_1}r(g(Y,Z)X&-g(X,Z)Y)\nnm\\
&-[\bar A^{\mathcal H^\bot}_0(X),\bar A^{\mathcal H^\bot}_0(Y)](Z)=0\label{barsigma0hbot-gauss}
\end{align}
for all $X,Y,Z\in {\mathcal H^\bot}$.

In deed, the fact that $A^0_{00}\in \ol{\mathcal S}_{\bbr^{q}}(L_1)$ implies
\be\label{td sigma0 gauss}
L_1(g(Y,Z)X-g(X,Z)Y)-[ A^0_{00}(X), A^0_{00}(Y)](Z)=0.
\ee

But by the decomposition \eqref{decomp tdsigma0}
\begin{align}
 A^0_{00}(&X)( A^0_{00}(Y)(Z)) = A^0_{00}(X, A^0_{00}(Y,Z))\nnm\\
=& A^0_{00}(X,\bar A^{\mathcal H^\bot}_0(Y,Z) +\bar A^{\mathcal H}_0(Y,Z))\nnm\\
=& A^0_{00}(X,\bar A^{\mathcal H^\bot}_0(Y,Z)) + A^0_{00}(X,\bar A^{\mathcal H}_0(Y,Z))\nnm\\
=&\bar A^{\mathcal H^\bot}_0(X,\bar A^{\mathcal H^\bot}_0(Y,Z)) +\bar A^{\mathcal H}_0(X,\bar A^{\mathcal H^\bot}_0(Y,Z))
+ A^{\mathcal H^\bot}_0(X,\bar A^{\mathcal H}_0(Y,Z))\nnm\\
=&\bar A^{\mathcal H^\bot}_0(X)(\bar A^{\mathcal H^\bot}_0(Y)(Z)) -\fr1rg(X,\bar A^{\mathcal H^\bot}_0(Y,Z))\sum_\alpha (n_\alpha+1)H_\alpha\nnm\\ &\quad-\fr1rg(Y,Z)\sum_\alpha (n_\alpha+1) A^{\mathcal H^\bot}_0(X,H_\alpha)\nnm\\
=&\bar A^{\mathcal H^\bot}_0(X)(\bar A^{\mathcal H^\bot}_0(Y)(Z))
-\fr1rg(X,\bar A^{\mathcal H^\bot}_0(Y,Z))\sum_\alpha (n_\alpha+1)H_\alpha\nnm\\ &\quad
-\fr1r\sum_\alpha (n_\alpha+1)L_1g(Y,Z)X\nnm\\
=&\bar A^{\mathcal H^\bot}_0(X)(\bar A^{\mathcal H^\bot}_0(Y)(Z))
-\fr1rg(X,\bar A^{\mathcal H^\bot}_0(Y,Z))\sum_\alpha (n_\alpha+1)H_\alpha\nnm\\ &\quad
-\fr{n-r+1}rL_1g(Y,Z)X\nnm
\end{align}
where we have used \eqref{000}, \eqref{barsigma0h-1} and the definition of $r$.

Since $g(X,\bar A^{\mathcal H^\bot}_0(Y,Z))=g(Y,\bar A^{\mathcal H^\bot}_0(X,Z))$, we find that
\begin{align}
[ A^0_{00}&(X), A^0_{00}(Y)](Z) = A^0_{00}(X)( A^0_{00}(Y)(Z)) - A^0_{00}(Y)( A^0_{00}(X)(Z))\nnm\\
=&\bar A^{\mathcal H^\bot}_0(X)(\bar A^{\mathcal H^\bot}_0(Y)(Z))
-\bar A^{\mathcal H^\bot}_0(Y)(\bar A^{\mathcal H^\bot}_0(X)(Z))\nnm\\
&\ -\fr{n-r+1}rL_1(g(Y,Z)X-g(X,Z)Y)\nnm\\
=&[\bar A^{\mathcal H^\bot}_0(X),\bar A^{\mathcal H^\bot}_0(Y)](Z)
-\fr{n-r+1}rL_1(g(Y,Z)X-g(X,Z)Y).
\end{align}
Inserting the above equality into \eqref{td sigma0 gauss} we obtain the equation \eqref{barsigma0hbot-gauss}, which completes the proof of Lemma \ref{barsigma0}.\endproof

\begin{lem}\label{added lem}
All the vectors $H_0,H_1,\cdots,H_s$ are constant vectors in ${\mathcal H}\subset\bbr^q$; In particular, $\bar c_0$ ,$\bar c_1$, $\cdots$, $\bar c_s$ are all constants. Furthermore,
\be\label{alp alp alp}
A^\alpha_{\alpha\alpha}\in {\mathcal S}_{ M_\alpha}(L_1-\bar c^2_\alpha).
\ee
\end{lem}

{\it Poof of Lemma \ref{added lem}}: The first conclusion of the lemma is easily obtained by taking the derivatives of \eqref{H0-1} on $M_\alpha$ for each $\alpha=0$, $1$, $\cdots$, $s$.

We first
note that the conditions (1) and (3) for \eqref{alp alp alp} come correspondingly from those for $A\in {\mathcal S}_{(M^n,g)}(L_1)$.

To verify the condition (2) for \eqref{alp alp alp}, we use the affine Gauss equation \eqref{gaus_af sph} to see that
\begin{align}
R^{M_\alpha}(X_\alpha,Y_\alpha)Z_\alpha =&L_1(g_{\alpha}(Y_\alpha,Z_\alpha)X_\alpha -g_{\alpha}(X_\alpha,Z_\alpha)Y_\alpha)\nnm\\
&\ \ \ \ -[A(X_\alpha),A(Y_\alpha)](Z_\alpha).
\label{added f4}
\end{align}
But by \eqref{decomposition of A}, \eqref{alp alp 0} and \eqref{alp 0 alp} we find
\begin{align*}
&[A(X_\alpha),A(Y_\alpha)](Z_\alpha)=A(X_\alpha,A(Y_\alpha,Z_\alpha)) -A(Y_\alpha,A(X_\alpha,Z_\alpha))\\ =&A(X_\alpha,A^0_{\alpha\alpha}(Y_\alpha,Z_\alpha) +A^\alpha_{\alpha\alpha}(Y_\alpha,Z_\alpha))\\ &-A(Y_\alpha,A^0_{\alpha\alpha}(X_\alpha,Z_\alpha) +A^\alpha_{\alpha\alpha}(X_\alpha,Z_\alpha)) \\
=&A^\alpha_{\alpha0}(X_\alpha,g_{\alpha}(Y_\alpha,Z_\alpha)H_\alpha) +A^0_{\alpha\alpha}(X_\alpha,A^\alpha_{\alpha\alpha}(Y_\alpha,Z_\alpha))\\ &+A^\alpha_{\alpha\alpha}(X_\alpha,A^\alpha_{\alpha\alpha} (Y_\alpha,Z_\alpha))
-A^\alpha_{\alpha0}(Y_\alpha,g_{\alpha}(X_\alpha,Z_\alpha)H_\alpha)\\ &-A^0_{\alpha\alpha}(X_\alpha,A^\alpha_{\alpha\alpha}(Y_\alpha,Z_\alpha))
-A^\alpha_{\alpha\alpha}(Y_\alpha,A^\alpha_{\alpha\alpha} (X_\alpha,Z_\alpha))\\
=&g_{\alpha}(Y_\alpha,Z_\alpha)g_\alpha(H_\alpha,H_\alpha)X_\alpha +g_\alpha(X_\alpha,A^\alpha_{\alpha\alpha}(Y_\alpha,Z_\alpha)H_\alpha\\ &+A^\alpha_{\alpha\alpha}(X_\alpha,A^\alpha_{\alpha\alpha} (Y_\alpha,Z_\alpha))
-g_{\alpha}(X_\alpha,Z_\alpha)g_\alpha(H_\alpha,H_\alpha)Y_\alpha\\ &-g_\alpha(Y_\alpha,A^\alpha_{\alpha\alpha}(X_\alpha,Z_\alpha)H_\alpha -A^\alpha_{\alpha\alpha}(Y_\alpha,A^\alpha_{\alpha\alpha} (X_\alpha,Z_\alpha))\\
=&\bar c^2_\alpha(g_{\alpha}(Y_\alpha,Z_\alpha)X_\alpha -g_{\alpha}(X_\alpha,Z_\alpha)Y_\alpha) +[A^\alpha_{\alpha\alpha}(X_\alpha),A^\alpha_{\alpha\alpha} (Y_\alpha)](Z_\alpha).
\end{align*}
Inserting the above into \eqref{added f4} we obtain
\begin{align*}
R^{M_\alpha}(X_\alpha,Y_\alpha)Z_\alpha =&(L_1-\bar c^2_\alpha)(g_{\alpha}(Y_\alpha,Z_\alpha)X_\alpha -g_{\alpha}(X_\alpha,Z_\alpha)Y_\alpha)\\
&-[A^\alpha_{\alpha\alpha}(X_\alpha),A^\alpha_{\alpha\alpha} (Y_\alpha)](Z_\alpha),
\end{align*}
implying the condition (2) for $A^\alpha_{\alpha\alpha}\in {\mathcal S}_{ M_\alpha}(L_1-\bar c^2_\alpha)$.\endproof

\begin{lem}\label{lem4.3}
The vector-valued symmetric bilinear form $ A\in {\mathcal S}_{(M^n,g)}(L_1)$ is uniquely, up to equivalence, determined by the metrics $ g_\alpha$, the flat metric $g_0$ on $\bbr^{q}$, the bilinear forms $ A^{\alpha}_{\alpha\alpha}$ $($$\alpha=1,\cdots,s$$)$ and the affine mean curvature $L_1$.\end{lem}

{\it Proof of Lemma \ref{lem4.3}}:
Since $ A^{\alpha}_{\alpha\alpha}\in {\mathcal S}_{( M_\alpha, g_\alpha)}(L_1-\bar c^2_\alpha)$, we see that the constant $\bar c_\alpha$ is completely determined by $g_\alpha$, $L_1$ and $ A^\alpha_{\alpha\alpha}$ via
\begin{align*}
R_{ g_\alpha}(X_\alpha,Y_\alpha)Z_\alpha=&(-L_1+ \bar c^2_\alpha)( g_\alpha(Y_\alpha,Z_\alpha)X_\alpha- g_\alpha(X_\alpha,Z_\alpha)Y_\alpha) \\ &+[ A^\alpha_{\alpha\alpha}(X_\alpha), A^\alpha_{\alpha\alpha}(Y_\alpha)]Z_\alpha,
\end{align*}
where $X_\alpha,Y_\alpha,Z_\alpha\in TM_\alpha$ and $R_{ g_\alpha}$ is the curvature tensor of $g_\alpha$.

On the other hand, up to an orthogonal transformation on ${\mathcal H}\subset TM_0\equiv\bbr^{q}$, the constant vectors $H_1,\cdots,H_s$ are uniquely given by the matrix equality
$$
\lmx g(H_1,H_1)&g(H_1,H_2)&\cdots&g(H_1,H_s)\\
g(H_2,H_1)&g(H_2,H_2)&\cdots&g(H_2,H_s)\\
\cdots&\cdots&\cdots&\cdots\\
g(H_s,H_1)&g(H_s,H_2)&\cdots&g(H_s,H_s)\rmx
=
\lmx  \bar c^2_1&L_1&\cdots&L_1\\
L_1& \bar c^2_2&\cdots&L_1\\
\cdots&\cdots&\cdots&\cdots\\
L_1&L_1&\cdots& \bar c^2_s\rmx.
$$

Furthermore, it is easily seen from \eqref{alp alp 0}, \eqref{alp 0 alp}, \eqref{000} and \eqref{barsigma0h} that
$ A^0_{\alpha\alpha}$, $ A^\alpha_{0\alpha}$, $ A^\alpha_{\alpha0}$, $\bar A^{\mathcal H}_0$, $\left. A_0^{\mathcal H^\bot}\right|_{{\mathcal H^\bot}\times {\mathcal H}}$, $\left. A_0^{\mathcal H^\bot}\right|_{{\mathcal H}\times {\mathcal H^\bot}}$ and $\left. A_0^{\mathcal H}\right|_{{\mathcal H}\times {\mathcal H}}$ are completely determined by the flat metric $g_0$, the vectors $H_1,\cdots,H_s$ and the affine mean curvature $L_1$.

Finally, since $\bar A^{\mathcal H^\bot}_0\in {\mathcal S}_{\bbr^{r-1}}\left(\fr{(n+1)L_1}r\right)$, it is realized as the Fubini-Pick form of a flat hyperbolic affine hypersphere in $\bbr^r$.
Then a theorem of L. Vranken, A-M. Li and U. Simon in \cite{vra-li-sim91} (also see \cite{amli89}) assures that any of such flat affine hypersphere is equiaffine equivalent to the hyperbolic hypersphere in Example \ref{expl}. Thus $\bar A^{\mathcal H^\bot}_0$ is also unique up to isometries on $\bbr^{r-1}$. It then follows that $ A^0_{00}$ is completely determined by the flat metric $g$, the sections $H_1$, $\cdots$, $H_s$ and the affine mean curvature $L_1$ up to isometries on $\bbr^{q}$.

Summing up, we have proved the conclusion of Lemma \ref{lem4.3}.\endproof

Now we are in a position to complete the proof of Theorem \ref{the main}.

Let $C$ be given by \eqref{newl1c}. Suitably choosing the constants $c_a (1\leq a\leq r)$, $c_{r+\alpha},\, \lalp (1\leq\alpha\leq s)$, we can also assume the first equality. For each $\alpha=1,\cdots,s$, fix one Riemannian metric $$ \galp=\fr{(n+1)L_1}{(n_\alpha+1)\lalp} g_\alpha$$ on $ M_\alpha$. Then by \eqref{alp alp alp}, have find that
$$ A^\alpha_{\alpha\alpha}\in {\mathcal S}_{( M_\alpha, \galp)}\left(\fr{(n_\alpha+1)(L_1-\bar c_\alpha^2)} {(n+1)L_1}\lalp\right).$$

We claim that
\be\label{claim}
\fr{(n_\alpha+1)(L_1-\bar c_\alpha^2)}{(n+1)L_1}=1,\mb{\ that is,\ } \bar c^2_\alpha=\fr{n-n_\alpha}{n_\alpha+1}(-L_1).
\ee

In fact, multiplying $h_{\alpha\gamma}$ to the both sides of \eqref{*} and then taking sum over $\alpha$ we have
\be\label{**}
1=\sum_{\alpha,\beta}h^{\alpha\beta}h_{\alpha\gamma} =-\sum_\alpha\fr{n_\alpha+1}{rL_1}h_{\alpha\gamma}.
\ee
Since, by \eqref{alp alp 0}, $h_{\gamma\gamma}= \bar c^2_\gamma$ and $h_{\alpha\gamma}=L_1$ for $\alpha\neq\gamma$, the right hand side of \eqref{**} is
\begin{align}
-\sum_\alpha\fr{n_\alpha+1}{rL_1}h_{\alpha\gamma}=&-\fr{n_\gamma+1}{rL_1} \bar c^2_\gamma-\sum_{\alpha\neq\gamma}\fr{n_\alpha+1}{rL_1}L_1\nnm\\
 =&-\fr{n_\gamma+1}{rL_1} \bar c^2_\gamma-\fr1r\sum_{\alpha\neq\gamma}(n_\alpha+1)\nnm\\
 =&-\fr{n_\gamma+1}{rL_1} \bar c^2_\gamma+\fr1r(n_\gamma+1)-\fr1r\sum_\alpha (n_\alpha+1)\nnm\\
 =&-\fr{n_\gamma+1}{rL_1} \bar c^2_\gamma+\fr1r(n_\gamma+1)-\fr1r(n-r+1)\nnm\\
 =&-\fr{n_\gamma+1}{rL_1} \bar c^2_\gamma-\fr1r(n-n_\gamma)+1\label{***}
\end{align}

From \eqref{**} and \eqref{***} we easily prove the claim \eqref{claim}.

The equality \eqref{claim} shows that $\Aalp\equiv A^\alpha_{\alpha\alpha}\in {\mathcal S}_{( M_\alpha, \galp)}(\lalp)$. It follows from Corollary \ref{cor2.1} that, for each $\alpha=1,\cdots,s$, there exists a hyperbolic affine hypersphere $x_\alpha:M^{n_\alpha}_\alpha\to \bbr^{n_\alpha+1}$ having $\stx{\alpha}{g}$, $\lalp$ and $\stx{\alpha}{A}$ as its affine metric, affine mean curvature and Fubini-Pick form respectively.

Suitably choosing the parameters $t^1,\cdots,t^{K-1}$, $K=r+s$, the original flat metric $g_0$ on $\bbr^{K-1}$ can be written as $g_0=\sum_{\lambda,\mu}g_{\lambda\mu} dt^\lambda dt^\mu$ where $g_{\lambda\mu}$ is defined by \eqref{g-lam mu}.

Now we consider the Calabi composition $\bar x$ of $r$ points and the $s$ hyperbolic affine hyperspheres $x_\alpha$, with the previously chosen constants $c_a,c_{r+\alpha}$ mentioned.  Then it follows that the original hyperbolic affine hypersphere $x$ is equiaffine equivalent to the Calabi composition $\bar x$ since they have the same affine metric and Fubini-Pick form by Proposition \ref{corr0}.
\endproof

As an application of Theorem \ref{the main}, we can easily recover the following result in a direct manner:

\begin{cor}\label{cor} (\cite{lix13})
A locally strongly convex and affine symmetric hypersurface $x:M^n\to\bbr^{n+1}$ is locally affine equivalent to the Calabi composition of some hyperbolic affine hyperspheres possibly including point factors if and only if $M^n$ is reducible as a Riemannian manifold with respect to the affine metric.\end{cor}

In fact, if $x$ is affine symmetric, then the holonomy algebra ${\mathfrak h}$ for $(M^n,g)$ acts trivially on the Fubini-Pick form $A$ which, together with the fact that, for each $\alpha=1,\cdots,s$, ${\mathfrak h}_\alpha$ acts irreducibly on $TM_\alpha$ and trivially on other $TM_\beta$ ($\beta\neq\alpha$), directly implies (2) and (3) in Theorem \ref{the main}.


\begin{thebibliography}{00}

\bibitem{bok-nom-sim90} N. Bokan, K. Nomizu and U. Simon, Affine hypersurfaces with parallel cubic forms, T\^ohoku Math. J. 42 (1990), 101-108, MR 1036477, Zbl0696.53006.
\bibitem{cal72} E. Calabi, {Complete affine hypersurfaces I}, Symposia Math., 10(1972), 19-38.
\bibitem{dil-vra94} F. Dillen and L. Vrancken, {\it Calabi-type composition of affine spheres}, Diff. Geom. appl, 4(1994), 303-328.
\bibitem{dil-vra-yap94} F. Dillen, L. Vrancken and S. Yaprak, Affine hypersurfaces with parallel cubic form, Nagoya Math. J. 135 (1994), 153-164. MR 1295822, Zbl0806.53008.306
\bibitem{hu-li-vra08} Z. J. Hu, H. Z. Li and L. Vrancken, {\it Characterizations of the Calabi product of hyperbolic affine hyperspheres}, Result. Math. 52 (2008), 299每314.
\bibitem{hu-li11} Z. J. Hu, C.C. Li, {\it The classification of $3$-dimensional Lorentian affine hypersurfaces with parallel cubic form}, Result. Math. 52(2008), 299每314.
\bibitem{hu-li-li-vra11a} Z. J. Hu, C.C. Li, H. Z. Li and L. Vrancken, {\it The classification of $4$-dimensional nondegenerate affine hypersurfaces with parallel cubic form}, Journal of Geometry and Physics, 61(2011), 2035每2057.
\bibitem{hu-li-li-vra11b} Z. J. Hu, C.C. Li, H. Z. Li and L. Vrancken, {\it Lorentzian affine hypersurfaces with parallel cubic form}. Res. Math., 59(2011), 577每620.
\bibitem{hu-li-sim-vra09}Z. Hu, H. Li, U. Simon, and L. Vrancken, {\it On locally strongly convex affine hypersurfaces with parallel cubic form, I}, Diff. Geom. Appl. 27(2009), no2, 188-205.
\bibitem{hu-li-vra11} Z. J. Hu, H. Li, L. Vrancken, {\it Locally strongly convex affine hypersurfaces with parallel cubic form}, J. Diff. Geom., 87(2011), 239-307.
\bibitem{kri-vra99} C. P. Wang, {\it Lorentian affine hyperspheres with constant affine sectional curvature}, Trans. Amer. Math. Soc., 352(1999), no4, 1581-1599.
\bibitem{amli89} A-M. Li, {\it Some theorems in affine differential geometry}, Acta Math. Sinica, N.S.,
5(1989), 345-354.
\bibitem{amli90} A-M. Li, {\it Calabi conjecture on hyperbolic affine hyperspheres}, Math. Z. 203(1990), 483-491.
\bibitem{amli92} A-M. Li, {\it Calabi conjecture on hyperbolic affine hyperspheres (2)}, Math. Ann. 293(1992), 485-493.
\bibitem{li-sim-zhao93} A-M. Li, U. Simon and G. S. Zhao, {\it Global affine differential geometry of hypersurfaces}, de Gruyter Expositions in Mathematics, vol. 11, Walter de Gruyter and Co., Berlin, 1993.
\bibitem{li-wang11} H.Z., Li and X.F. Wang, {\it Calabi product Lagrangian immersions in complex projective space and complex hyperbolic space}. Results Math., 59(2011), 453-470.
\bibitem{lix93} X. X. Li, {\it The composition and the section of hyperbolic affine spheres}, J. Henan Normal University (Natural Science Edition, in Chines), 21(1993), no.2, 8-12.
\bibitem{lix11} X. X. Li, {\it On the Calabi composition of multiple affine hyperspheres}, preprint, 2011.
\bibitem{lix13} X. X. Li, {\it On the correspondence between symmetric equiaffine hyperspheres
and the minimal symmetric Lagrangian submanifolds},  preprint in Chinese, 2013. To appear.
\bibitem{nom-sas94} K. Nomizu and T. Sasaki, Affine Differential Geometry. Cambridge University Press, Cambridge (1994).
\bibitem{vra-li-sim91} L. Vrancken, A-M. Li and U. Simon, {\it Affine spheres with constant affine sectional curvature}, Math. Z. 206(1991), 651-658.
\bibitem{wang93} C. P. Wang, {\it Canonical equiaffine hypersurfaces in $\bbr^{n+1}$}, Math. Z., 214(1993), 579-592.

\end{thebibliography}
\end{document}